\documentclass{amsart}
\usepackage{amsfonts,amscd,amssymb,amsmath,amsthm}
\usepackage{graphicx,mathrsfs,appendix}
\usepackage{centernot,caption,subcaption,color,url}
\usepackage[T1]{fontenc}
\usepackage{ctable}
\usepackage[a4paper, margin=1in]{geometry}

\usepackage{enumitem}
\usepackage[colorlinks=true,citecolor=blue]{hyperref}
\usepackage{cleveref}
\usepackage{cite}

\begin{document}

\newtheorem{theorem}{Theorem}
\newtheorem{lemma}[theorem]{Lemma}
\newtheorem{proposition}[theorem]{Proposition}
\newcommand{\R}{\mathbb{R}}
\newcommand{\N}{\mathbb{N}}
\newcommand{\T}{\mathcal{T}}
\newcommand{\B}{\mathcal{B}}
\newcommand{\A}{\mathcal{A}}
\newcommand{\G}{\mathcal{G}}
\newcommand{\F}{\mathcal{F}}
\newcommand{\Z}{\mathbb{Z}}
\newcommand{\Q}{\mathbb{Q}}
\newcommand{\1}{\mathbf 1}

\newcommand{\limd}[1]{\underset{#1}{\lim }\,}
\renewcommand{\Pr}[2][{}]{\mathbb{P}_{#1}\left(#2\right)}
\newcommand{\limn}[0]{\underset{n\to\infty}{\lim}\,}
\newcommand{\E}[2][{}]{\mathbb{E}_{#1}\hspace{-0.15em}\left[#2\right]}
\newcommand{\sumdu}[2]{\overset{#2}{\underset{#1}{\sum}}\,}
\newcommand{\proddu}[2]{\overset{#2}{\underset{#1}{\prod}}}
\newcommand{\intdu}[2]{\int_{#1}^{#2}}
\newcommand\independent{\protect\mathpalette{\protect\independenT}{\perp}}
\def\independenT#1#2{\mathrel{\rlap{$#1#2$}\mkern2mu{#1#2}}}

\title[Inferring cumulative advantage from longitudinal records]{Inferring cumulative advantage from longitudinal records}
\author{Alexandros Gelastopoulos}
\address{Alexandros Gelastopoulos, University of Southern Denmark, Pompeu Fabra University, Barcelona School of Economics}
\email{alex@sam.sdu.dk}
\author{Lucas Sage}
\address{Lucas Sage, European University Institute}
\email{lucaschristophemarc.sage@eui.eu}
\author{Arnout van de Rijt}
\address{Arnout van de Rijt, European University Institute}
\email{arnout.vanderijt@eui.eu}

\thanks{We'd like to thank Damon Centola, Eran Shor, and Mathijs de Vaan for helpful comments on an earlier draft, and members of the Colloquium on Analytical Sociology at the European University Institute. A. Gelastopoulos was supported in part from Sapere Aude grant 2065-00038B by the Independent Research Fund Denmark, and from ERC Consolidator Grant \#772268 from the European Commission.}
\date{\today}

\begin{abstract}
Inequality in human success may emerge through endogenous success-breeds-success dynamics but may also originate in pre-existing differences in talent. It is widely recognized that the skew in static frequency distributions of success implied by a cumulative advantage model is also consistent with a talent model. Studies have turned to longitudinal records of success, seeking to exploit the time dimension for adjudication. Here we show that success histories suffer from a similar identification problem as static distributional evidence. We prove that for any talent model there exists an analogous path dependent model that generates the same longitudinal predictions, and vice versa. We formally identify such twins for prominent models in the literature, in both directions. These results imply that longitudinal data previously interpreted to support a talent model equally well fits a model of cumulative advantage and vice versa.
\end{abstract}

\maketitle

\newcommand{\comment}[1]{\textcolor{red}{(#1)}}

\bibliographystyle{unsrt}

A general hypothesis across the sciences is that at the origins of many observed inequalities lies a generative process of positive reinforcement, a.k.a. cumulative advantage, the Matthew Effect, success-breeds-success or rich-get-richer \cite{merton1968matthew, price1976general, diprete2006cumulative}. In theoretical models, cumulative advantage is capable of producing persistent success differentiation between ex ante equivalent individuals, through path dependence on early random success. The prevalence of positive feedback in success accumulation in many contexts is made plausible by the existence of a number of generic mechanisms that could each generate it, such as status \cite{correll2017s}, diffusion \cite{young2009innovation}, social learning \cite{bikhchandani1992theory, salganik2006experimental}, and increasing returns \cite{arthur1989competing, dimaggio2012network}. A well-known inference problem is that the distributional skew implied by models of cumulative advantage that is robustly found across contexts is also consistent with competing generative models that lack a reinforcement element, e.g. a model that assigns each individual an own propensity for success (a `talent') \cite{coleman1964introduction, albert2002statistical, newman2005power, simon1955class, mandelbrot1959note, mitzenmacher2004brief, stumpf2012critical, perc2014matthew, denrell2012top}. Static distributions are therefore a weak form of evidence for path dependence in success accumulation.

Faced with this inference problem, it is natural to turn to longitudinal data and seek to exploit the time dimension for adjudication between cumulative advantage and talent as generative principles: Cumulative advantage, and more generally path dependence, links outcomes across time, while talent does not. A large body of literature thus investigates longitudinal patterns of earnings, learning outcomes, publications and citations, interpreting diverging within- or between-group trends as evidence for cumulative advantage and time-constant differences as consistent with a talent model \cite{allison1982cumulative, diprete2006cumulative, wagner1995were, huber1998cumulative, sinatra2016quantifying}.

In this paper we find that dynamic records of success suffer from a similar identifiability problem as static data. Others have previously demonstrated  issues with differentiating between the two models in dynamic data \cite{allison1980estimation, cane1977class,lancaster1980analysis,elbers1982true,heckman1984identifiability}. We generalize prior results to prove that for any talent model there exists an analogous path dependent model that generates the same longitudinal predictions. Conversely, we prove that for any path dependent model there exists an analogous talent model generating the same longitudinal predictions. We formally identify the cumulative advantage equivalent of the $Q$-model of citation counts \cite{sinatra2016quantifying}, and the talent model equivalent of the contagious Poisson process \cite{allison1982cumulative}. These results imply that all existing observational studies reporting evidence in favor of path dependence must equally support a talent model, and vice versa.

\section*{Theoretical results}

\subsection*{Cumulative advantage vs talent in binary outcomes}

A standard model of cumulative advantage is the P\'olya urn model, which can be described as follows: An urn contains balls of various colors. At each time step, we randomly draw one ball from the urn and then return it together with another ball of the same color. That is, we add one ball at each time step, with the probability of that ball having color $i$ proportional to the number of balls of that color $i$ already in the urn. In the most basic version of the model, no new colors can appear. Modifications of this scheme, often allowing the appearance of new colors, form the basis of many other cumulative advantage models, such as Simon's reinforcement model \cite{simon1955class}, Barab\'asi and Albert's preferential attachment model \cite{barabasi1999emergence}, and Price's urn model \cite{price1976general} (see the Appendix for a discussion of these models).

If the urn initially contains balls of two colors only, say one black and one red ball, then this is a simple model of success-breeds-success. The draws represent a series of failures and successes in the endeavours of an individual, such as publications of a researcher, works of an artist, etc., with black balls interpreted as successes and red balls as failures. Given that the more often we observe a success the more likely successes become in the future, this is a cumulative advantage model. However, it turns out that this model is statistically equivalent with a pure talent model.

Specifically, consider the following alternative model. An individual has an inherent talent $p$, which is a number drawn from a given distribution on the interval $[0,1]$, and it remains fixed throughout their lifetime. $p$ is interpreted as the probability of success of each of the individual's outputs (e.g. a paper for a researcher). More precisely, we assume that we observe a series of binary draws, successes or failures, with the probability of each draw being a success equal to $p$, independently of the outcomes of previous draws\footnote{Here we refer to \textit{causal} independence. The draws are not statistically independent, because they are linked through the value of the random variable $p$. See next section.}.

Assuming that $p$ is drawn uniformly randomly from the interval $[0,1]$, and that it cannot be \textit{directly} observed, then by observing a sequence of failures and successes, it is impossible to infer whether this sequence was generated by a P\'olya urn model or by the above talent model. This is formalized in the following proposition, which draws on De Finetti's theorem \cite{de1931funzione,de1937prevision} (see for example \cite[Theorems 2.1 and 2.2]{freedman1965bernard} or \cite[Example 1.5]{mahmoud2008polya}).
\begin{proposition}
\label{propositionPolyaUrn}
    Consider a P\'olya urn model starting with 1 black and 1 red ball, and let $Y_n$ denote the outcome of the $n$-th draw with $Y_n=1$ if a black ball is drawn and $Y_n=0$ otherwise. Also, let $p$ be uniformly distributed on $[0,1]$ and let $Y_n'$ be a sequence of binary random variables, conditionally independent given $p$, with $\Pr{Y_n'=1\mid p}=p$. Then, the sequences $\{Y_1,Y_2,\ldots \}$ and $\{Y'_1,Y'_2,\ldots \}$ have the same joint probability distribution, i.e. they are statistically indistinguishable.
\end{proposition}
We emphasize that the proposition pertains to \textit{joint} probability distributions, the distributions of entire sequences, not of single observations. Equality of joint distributions is a much stronger assertion than equality of the distributions of $Y_n$ and $Y_n'$ for each $n$ separately. Two sequences having identical probability distributions implies that there is no statistical test that can decide which of the two generative models produced a given empirical sequence or even say which is more likely. This claim continues to hold even if we observe series of successes and failures of multiple individuals;\footnote{Assuming, for example, that the sequences of successes for different individuals are independent.} it is impossible to say whether these series of observations were all generated by multiple realizations of the P\'olya urn model (which contains no notion of talent) or by a talent model in which each individual has an unobserved talent $p_i$, drawn uniformly from $[0,1]$. This counterintuitive result tells us that the distinction between a talent model and a cumulative advantage process in this case is merely epistemic, it cannot be detected in the data.

\subsection*{Path dependent equivalents of talent models}

We now argue that the above equivalence between a talent model and a cumulative advantage model is not an exceptional case, nor does it hold only for processes with binary observations. We will show that for a general class of generative models that rely on an unobserved talent variable, one can construct a path dependent process\footnote{We will define a path dependent model as one in which future performance depends causally on past performance. Cumulative advantage is a special case of path dependence where each success makes future successes more likely. To the best of our knowledge, there is no precise mathematical definition of cumulative advantage in a general setting. Since our results are stated in terms of the more general notion of path dependence, we do not provide such a definition either. As will become evident, the examples we give can naturally be interpreted as cumulative advantage models.} that contains no notion of talent, yet produces the exact same statistics. A partial converse also holds.

We first consider a case of a single individual observed at just two time points, which we then generalize. Let $T,X_1,X_2$ be three independent standard normal random variables, i.e. with zero mean and unit variance. Suppose that we do not directly observe their values, but instead we observe $Y_1=T+X_1$ and $Y_2=T+X_2$. Here $T$ can be interpreted as the inherent talent of an individual in some skill, and $Y_1$ and $Y_2$ as the quality of their first and second work, respectively, which are affected by the chance components $X_1$ and $X_2$. Since $X_1$ and $X_2$ are independent, this model contains no cumulative advantage; $Y_2$ depends only on the individual's talent and their luck at time $t=2$; it is not causally affected by the degree of their success at time $t=1$.

Despite the fact that $Y_1$ and $Y_2$ are not causally related to each other, these variables are not independent in the statistical sense. Knowing the value of $Y_1$ provides information about the value of $T$, which in turn provides information about $Y_2$. Specifically, a larger $Y_1$ will imply a larger $T$ on average, which in turn implies a larger $Y_2$. $Y_1$ and $Y_2$ are positively correlated.

We now define another pair of random variables, $Y'_1$ and $Y'_2$, with the same joint distribution as $Y_1, Y_2$, but which exhibit cumulative advantage in the sense that $Y'_2$ is causally related to $Y'_1$. Specifically, let $X_1,X_2$ be as before, but now define $Y'_1=2X_1$ and $Y'_2=\frac{Y'_1}{2}+X_2$. This is a cumulative advantage model, because the performance at time $t=2$ is directly affected by the performance at $t=1$, with an increase in $Y'_1$ resulting in an increase in $Y'_2$. Additionally, this model contains no notion of talent. Despite their very different causal structures, the two pairs of random variables have the same statistics: for their second moments we have $Var(Y_1)=Var(Y_2)=2$, $Cov(Y_1,Y_2)=1$, and the same equalities hold for the pair $Y'_1,Y'_2$. Since both pairs are also jointly normal with zero means, their joint distributions are identical. Someone who observes only occurrences of $Y_1$ and $Y_2$ would have no way to tell whether they were produced by the first or the second model. The situation wouldn't change even if we observed multiple independent pairs of $Y$'s, i.e. $(Y^1_1,Y^1_2)$, $(Y^2_1,Y^2_2)$, \ldots , $(Y^n_1,Y^n_2)$; it would still be impossible to tell whether they were all produced by the first or all by the second model.

The above observations can be generalized. Let us define a \textit{talent model} as a procedure that generates a sequence of observations $\{Y_n\}$,\footnote{We use the notation $\{Y_n\}$ to refer to the entire sequence of observations, and $Y_n$ to refer to the $n$-th element of that sequence.} typically taking values in some subset of $\R$ (but not necessarily), whose $n$-th term is a function of (i) a ``talent'' variable $T$ and (ii) a chance component $X_n$, with the sequence $\{X_n\}$ being i.i.d. random variables. In order to account for multiple individuals, we will also assume that $T$ is a random variable that is drawn from a fixed distribution. In principle we make no assumptions about the form of the dependence of $Y_n$ on $T$ and $X_n$; that is, we let the $n$-th observation be given by $Y_n=f(T,X_n)$, where $f$ can be any function.\footnote{Technically, we require that $f$ be \textit{measurable}, which is the minimal requirement for $Y_n$ to be a random variable.} For example, in the previous example we had $f(T,X_n)=T+X_n$.

Next we describe what we mean by a purely path dependent model, of which a cumulative advantage model is a special case. A \textit{purely path dependent model} is one in which there is no notion of talent, but instead the $n$-th observation $Y_n$ depends on a combination of luck (modelled as an i.i.d. sequence $\{X_n\}$, as before) and the values of the previous observations, $Y_1,\ldots ,Y_{n-1}$. Again, we make no assumptions about the form of these dependencies. That is,
\begin{equation}
\label{eqPathDependentProcess}
    Y_n=g_n(Y_1,\ldots ,Y_{n-1},X_n),
\end{equation}
where $g_n$ can be any measurable function. In the previous example we had $g_1(X_1)=X_1$ and $g_2(Y_1,X_2)=Y_1\slash 2+X_2$. In the Appendix we show how several well-known models of cumulative advantage, including preferential attachment models, multiplicative processes (i.e., processes that obey Gibrat's law), and Price's urn model can be expressed in terms of \cref{eqPathDependentProcess}.

The next theorem states our first main result, namely that for any talent model, there is a purely path dependent model that is statistically indistinguishable from the first.

\begin{theorem}
\label{theoremTalentToPathDependent}
For any talent model that generates a sequence $\{Y_n\}$, there is a purely path dependent model generating a sequence $\{Y'_n\}$ with the same joint distribution as $\{Y_n\}$.
\end{theorem}

The proof of this and all following theorems are given in the Appendix.

\subsection*{A concrete example: a purely path dependent equivalent of the $Q$-model of citations}

Sinatra et al. \cite{sinatra2016quantifying} propose the following talent model of citation counts: The number of citations $c_n$ of the $n$-th paper published by an individual is the product of a factor $Q$, which is fixed throughout the lifetime of the individual, and a luck component $p_n$, with the $p_n$'s forming an i.i.d. sequence.\footnote{Here we assume that infinitely many papers are produced by each individual, i.e. we ignore the question of productivity.} The random variable $Q$ and the $p_n$'s are assumed to be independent and each is log-normally distributed. As a result, $c_n=Q\cdot p_n$ is also log-normally distributed. For convenience we are going to work with the logarithms of these quantities, which are all normally distributed, and for which the relation $\log c_n=\log Q+\log p_n$ holds. This is a talent model as defined in the previous section, if we let $Y_n=\log c_n$, $T=\log Q$, $X_n=\log p_n$, and the function $f$ is taken to be $f(T,X_n)=T+X_n$, so that $Y_n=f(T,X_n)$. In \cite{sinatra2016quantifying}, the mean of $X_n$ is assumed to be $0$, and we follow the same convention.

For this talent model, the cumulative advantage (purely path dependent) model that produces the same statistics can be given explicitly as follows. Let $a\in\R$ and $b,c>0$ be some parameters, and define $Y_n'$ recursively so that conditioned on $Y_1',\ldots ,Y_n'$, the variable $Y_{n+1}'$ is normally distributed, with mean
\begin{equation}
\label{eqMean1}
    \mu_{n+1} = \frac{c\cdot  a+\sumdu{k=1}{n}Y_k'}{n+c}
\end{equation}
and variance
\begin{equation}
\label{eqVariance1}
    \sigma^2_{n+1} = b \cdot\left(1+\frac{1}{n+c}\right).\footnote{In terms of the formulation of the previous section, we may choose 
$$
g_{n+1}(Y_1',\ldots,Y_{n}',X'_{n+1})=\mu_{n+1}+\sigma_{n+1}\cdot X'_{n+1},
$$
where $X'_{n+1}$ is a standard normal random variable.
}
\end{equation}

In particular, for $n=0$ we get that the unconditional distribution of $Y'_1$ is ${\mathcal N}(a,b+b/c)$.

Note that in the above model there is no notion of talent. Instead, the value of $Y'_{n+1}$ is determined by the values of $Y'_1,\ldots ,Y_n'$ and noise. Specifically, $Y'_{n+1}$ is normally distributed, with known variance, and a mean that is a deterministic function of the past performances $Y_1',\ldots ,Y_n'$.

The next proposition states that, for matching values of parameters, this generative model is statistically indistinguishable from the $Q$-model of citations, i.e. any data that is consistent with one is also consistent with the other (see the Discussion for why one might prefer one model over the other nevertheless).

\begin{proposition}
\label{propositionQmodel}
    Let $\sigma_X$ and $\sigma_T$ be the standard deviation of $X_n$ and $T$ respectively, in the $Q$-model, and let $\mu_T$ be the mean of $T$. If we take $a=\mu _T$, $b=\sigma ^2_X$, and $c=\sigma^2_X/\sigma_T^2$, then the sequences $\{Y_n\}$ and $\{Y_n'\}$ have the same joint distribution.
\end{proposition}

\subsection*{Talent equivalents of  path dependent models}

We have established that every talent model can be mimicked by a purely path dependent model. We now show that the converse is also true if either we allow the function $f$ that maps talent and luck to output to vary over time or the path dependent process is \textit{exchangeable}.

Exchangeable means that the statistics of the process remain unchanged when the order of observations is shuffled: the sequence $Y_1,Y_2,\ldots $ is exchangeable if it is statistically indistinguishable from $Y_{\sigma _1},Y_{\sigma _2},\ldots $ for any finite permutation $\sigma_1,\sigma _2,\ldots $ of the positive integers.\footnote{A permutation is finite if it leaves all indices fixed except a finite number of them.} Every i.i.d. sequence is exchangeable, but an exchangeable sequence need not consist of independent random variables in general. An example of an exchangeable process that is not i.i.d. is the sequence of draws from a P\'olya urn \cite[Section 1.7]{mahmoud2008polya}. We have the following theorem.

\begin{theorem}
\label{theoremPathDependentToTalentExchangeable}
    Let $\{Y_n\}$ be any stochastic process. Then, there exists a talent model $Y'_n=f(T,X_n)$ such that $\{Y'_n\}$ has the same joint distribution as $\{Y_n\}$ if and only if $\{Y_n\}$ is exchangeable.
\end{theorem}

The reason that no talent model can mimic non-exchangeable path dependent models is that, as we have defined it above, the statistics of a talent model remain necessarily constant with $n$. Thus, if a stochastic process $\{Y_n\}$ has the property, for example, that $\E{Y_n}\neq \E{Y_m}$ for some $m\neq n$, then it is impossible for a talent model to statistically reproduce this process.

To allow for more flexibility, we define a \textit{time-varying talent model} as one in which $f$ can vary with time. That is, there is some talent variable $T$, some i.i.d. sequence $X_n$, and a sequence of functions $f_n$, so that $Y_n=f_n(T,X_n)$. With this more general definition, there is always a (time-varying) talent model that can mimic a path dependent model, and in fact any stochastic process.

\begin{theorem}
\label{theoremPathDependentToTalent}
    For any purely path dependent model that generates a sequence $\{Y_n\}$, there exists a time-varying talent model generating a sequence $\{Y'_n\}$ that has the same joint distribution as $\{Y_n\}$.
\end{theorem}

Although here we are focusing on the equivalence between purely path dependent and purely talent-based models, the assertion of \cref{theoremPathDependentToTalent} continues to hold if our starting model also incorporated talent, i.e. if the function $g_n$ in \cref{eqPathDependentProcess} also depended on $T$. Analogously, \cref{theoremTalentToPathDependent} continuous to hold if the starting model was time-varying and/or incorporated path dependence. The proofs would remain the same.\footnote{Conversely, given any purely path dependent or purely talent-based model, one can get any number of equivalent models with varying degrees of reliance on path dependence and talent. A simple way to do this is to first obtain two equivalent models, one of which is purely path dependent and the other purely talent-based (which we know exist by \cref{theoremTalentToPathDependent,theoremPathDependentToTalent}), and generate a sequence from the first one with probability $q$ and from the second otherwise.}

The proof of \cref{theoremPathDependentToTalent} shows that even a \textit{deterministic} talent model (i.e. one without the luck factor $X_n$) can produce any given distribution of observations. All randomness, correlations between observations, etc., can be attributed to an unobserved variable $T$. We note however that, although we have been calling $T$ ``talent'', there is no guarantee that this variable enjoys properties that we would normally associate with talent, e.g. that for given $Y_1,\ldots ,Y_{n-1}$, the value of $Y_n$ would be positively correlated with $T$. The theorem merely states that, as long as there are unobserved variables with unrestricted effect on the outcome, any statistical observations can be entirely attributed to those unobserved variables.

In practice, nonetheless, when given a specific path dependent model, we might be able to derive an explicit talent model equivalent to the path dependent one, where $T$ does behave as a ``talent'' variable, as shown next.

\subsection*{Talent equivalent of the contagious Poisson process model}

A popular model of ``event statistics'' is the contagious Poisson process \cite{coleman1964introduction,allison1982cumulative}, introduced by Greenwood and Yule \cite{greenwood1920inquiry}. Depending on the application, events may refer to illnesses \cite{eaton1974mental}, accidents \cite{arbous1951accident}, episodes of violence \cite{spilerman1970causes}, publications or citations \cite{allison1980estimation}. The term ``contagious'' refers to the fact that, unlike the standard Poisson process, the occurrence of each event makes future events more likely to occur.

More precisely, the contagious Poisson process is a continuous-time Markov process $X_t$ that measures the number of events that have occurred up to time $t$. The rate of new events at time $t$ is equal to $\alpha +\beta X_t$, for some positive constants $\alpha$ and $\beta$. Because the rate of events increases with the number of past events, this is a cumulative advantage process. Given that $\alpha $ and $\beta $ are constants (rather than random variables), the process is a \textit{purely} path dependent one. Because it is a continuous-time process, it cannot be expressed in the form of \cref{eqPathDependentProcess}, but the analogy should be clear.

It has long been known \cite{feller1943general} that by using static frequency data, the above model cannot be distinguished from a homogeneous mixed Poisson process \cite{grandell1997mixed,kallenberg2017random}, i.e. one that assumes a \textit{constant} rate of new events, but different across individuals, which is an instance of a (continuous-time) stationary talent model. Many researchers have thus suggested the use of longitudinal data and proposed a number of statistical methods to distinguish the contagious from the homogeneous mixed Poisson model \cite{coleman1964introduction,spilerman1970causes,feller1943general,eaton1978third,wasserman1983distinguishing,karlis1998estimation}. However, Lundberg \cite{lundberg1940random} showed that a time-discounted version of the contagious Poisson process, originally proposed by Eggenberger and P\'olya \cite{eggenberger1923statistik}, is indistinguishable from the mixed Poisson process \textit{also in analysis of longitudinal data}.\footnote{Cane \cite{cane1974concept} gave an independent proof many years later, unaware of Lundberg's work. See also \cite{taibleson1974distinguishing}.} Here we show that by ``inverting'' the time-discounting and applying it to the mixed Poisson process instead of the contagious one, we get that the original contagious process (without time-discounting) is indistinguishable from a time-varying talent model. Specifically, let $Q$ be a random variable distributed as $\text{Gamma}\left(\frac{\alpha}{\beta },\beta \right)$.\footnote{We use the shape-scale parameterization of the Gamma distribution, i.e. $Gamma(l,\theta)\sim \frac{1}{\Gamma(l)\theta^l}x^{l-1}e^{-x/\theta}$.} We define $N_t$ to be the number of events up to time $t$ of a (non-homogeneous) mixed Poisson process with rate given by $Q\cdot e^{\beta t}$. Interpreting $Q$ as talent,\footnote{We are using $Q$ instead of $T$ here, to avoid confusion with the time variable.} this model says that for a given value of talent ($Q=q$), the rate of new events is $r_q(t)=q\cdot e^{\beta t}$, i.e. it increases with time, but in a deterministic way; it does not depend on past events. Therefore, this is a continuous-time analogue of a time-varying talent model. We have the following proposition.

\begin{proposition}
\label{propAllisonTalentEquivalent}
    Let $X_t$ be a contagious Poisson process with parameters $\alpha,\beta >0$, and $N_t$ the number of events up to time $t$ of a (non-homogeneous) mixed Poisson process with rate given by $Q\cdot e^{\beta t}$, where $Q$ is $\text{Gamma}(\frac{\alpha}{\beta },\beta )$-distributed. Then, $\{N_t\}_{t\geq 0}$ and $\{X_t\}_{t\geq 0}$ have the same joint distribution.
\end{proposition}

As before, this proposition implies that by observing sequences of events, possibly from a cohort of individuals/units of observation, it is impossible to tell whether the variation is due to a contagious process or due to unobserved inherent heterogeneity combined with time-varying rates that are independent of past events (e.g., changing environment or career progress). Note that in certain settings it is possible that there is \textit{observable} heterogeneity, and the data may even exclude both further (unobservable) heterogeneity \textit{and} path dependence \cite{ritterband1973group}. But excluding only one of the two is impossible without further assumptions.\footnote{For example, Spilerman \cite{spilerman1970causes} distinguishes between heterogeneity and contagion based on the premise that either reinforcement or a non-path-dependent increase in the rate of incidents may take place, but not both at the same time.}

\subsection*{Path dependence across individuals}

The definition of a path dependent model we have given might at first sight seem to imply that the performance of an individual must depend only on the past performance of the \textit{same} individual (and luck). This is the case if $Y_n$ is taken to denote the performance of a single individual; then, our definition of a path dependent model as $Y_n=g_n(Y_1,\ldots ,Y_{n-1},X_n)$ implies that the performance of other individuals in a cohort cannot affect $Y_n$.

A simple model that does not satisfy this assumption is the following adaptation of a P\'olya urn: There are balls of three colors. When a ball of color 1 or 2 is drawn, we return it together with another ball of that color, but when color 3 is drawn, we return it with two extra balls of that color. In that case, the probability of drawing color 1 depends not only on how often color 1 has been drawn in the past, but also on how often each of the other colors has been drawn, because this will determine the total number of balls in the urn. If $Y_n$ is a binary variable that expresses whether a ball of color $1$ is drawn at time $n$, then a function of the form $g_n(Y_1,\ldots ,Y_{n-1},X_n)$ cannot express information about the draws of other colors. However, $Y_n$ can be taken instead to be a binary \textit{vector} whose $i$-th coordinate tells us whether a ball of color $i$ has been drawn. In that case, the sequence $Y_1,\ldots ,Y_{n-1}$ contains all information about previous draws, and consequently the value of $Y_n$ can be expressed in the form $g_n(Y_1,\ldots ,Y_{n-1},X_n)$.

More generally, in the definition of a path dependent model, $Y_n$ can be taken to contain information about the performance of all individuals in the cohort. Since \cref{theoremPathDependentToTalent} does not require $Y_n$ to have a specific form (e.g. to be real-valued), its assertion continuous to hold. In that case, however, the random variable $T$ given by the theorem is even less straightforward to interpret, since it contains information about all correlations between all individuals' performances at all times.

\section*{Empirical results}

To illustrate our theoretical findings, we use data from the MusicLab experiments \cite{salganik2006experimental}. Participants were allocated to nine distinct experimental conditions, referred to as `worlds', and invited to listen to and, if so desired, download one or more of 48 unfamiliar rock songs. In eight of the nine worlds, participants saw counts of prior downloads by earlier participants. In these influence worlds, both a talent model and a path dependent model are in principle appropriate, as it is mechanically possible that download behavior is largely predetermined by song quality and it is also possible that it is mainly driven by the idiosyncratic downloads of first movers and the cumulative advantage dynamics following that early history. Our focus instead will be on the ninth world, wherein individuals could not see what others had previously downloaded, so had no information about the opinions of their peers regarding song quality. In this independent condition, then, a path dependent model is unambiguously wrong, as a causal relationship between the download behavior of earlier and later subjects is precluded by design. Only a song's intrinsic appeal, i.e. its `talent', and its location on the screen, which was randomized, can determine song success. This then presents an ideal scenario for illustrating the inference problem identified in this article: Theorem 2 implies that if we can fit a talent model to this data, then it must also be possible to achieve as good a fit with a path dependent model. We aim to showcase that the data in fact permit an excellent fit to a simple P\'olya urn model, a model that we know is wrong.

The P\'olya urn model is implemented as follows: Each song corresponds to a unique ball color. When a user downloads a song, it corresponds to drawing a ball of the corresponding color. Initially, all colors have an equal probability of being drawn – the urn contains one ball of each color. Each ball drawn is put back into the urn, accompanied by \textit{f} other balls of the same color. We follow the empirical sequence of participant downloads. Each simulated user downloads as many songs as in the experiment. Because in the experiment participants are precluded from downloading the same song multiple times, also in our simulations, when a subject downloads more than one song, songs previously downloaded by that subject are assigned a zero probability, and other songs' probabilities proportionately inflated. (Simulation details are in the Appendix.)

\Cref{fig:musicLab1} charts the evolution of the percentage of downloads per song ranked by their popularity over time for the first version of the MusicLab experiment, in which songs were shown in a random position in a 16 x 3 rectangular grid. \Cref{fig:musicLab2} displays analogous results for the second version of the experiment, in which songs were presented in a single column in random order. The panels represent four snapshots: The first panel shows the download percentages during the first quarter of downloads in descending order, the second panel the second quarter, and so on. Blue dots are observed frequencies and green lines correspond to symmetric 95\% intervals in our simulations of a (wrong) P\'olya urn model.\footnote{The slight variations in the locations of the 95\% intervals stem from empirical differences across users in numbers of downloads.} \Cref{fig:musicLab1,fig:musicLab2} show that our wrong model is capable of generating a distribution of download percentages that provides an excellent fit to the empirical data. 

These results demonstrate that distinguishing between data generated by a cumulative advantage model and those arising from a talent model proves elusive. It is only because of the unique design of the MusicLab study, which precludes path dependence in the independent condition of the experiment, that we know that the cumulative advantage model that fits so well is wrong. Had the data come from an observational study we might have drawn the wrong conclusion about the origins of success in cultural markets.

\begin{figure*}[htp]
    \centering
    \includegraphics[width=15cm]{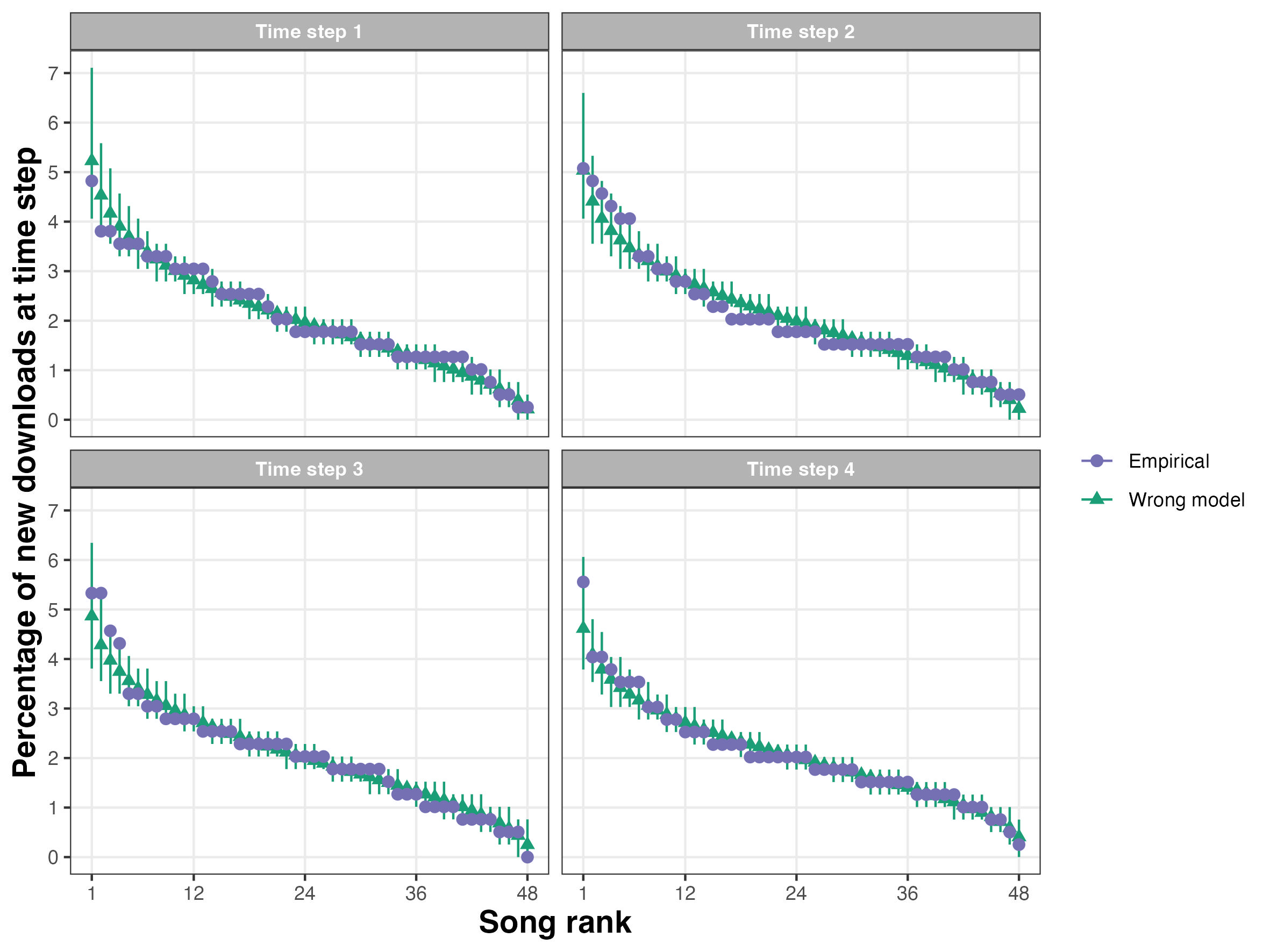}
    \caption{Downloads (\%) of 48 songs in MusicLab experiment version 1, ranked in descending order, by time step. Each time step contains a quarter of downloads in chronological order. Blue dots represent empirically observed percentages. Green triangles represent estimated percentages and green lines 95\% symmetric intervals from 2000 simulation runs of a P\'olya urn model with \textit{f}=0.295.}
    \label{fig:musicLab1}
\end{figure*}

\begin{figure*}[htp]
    \centering
    \includegraphics[width=15cm]{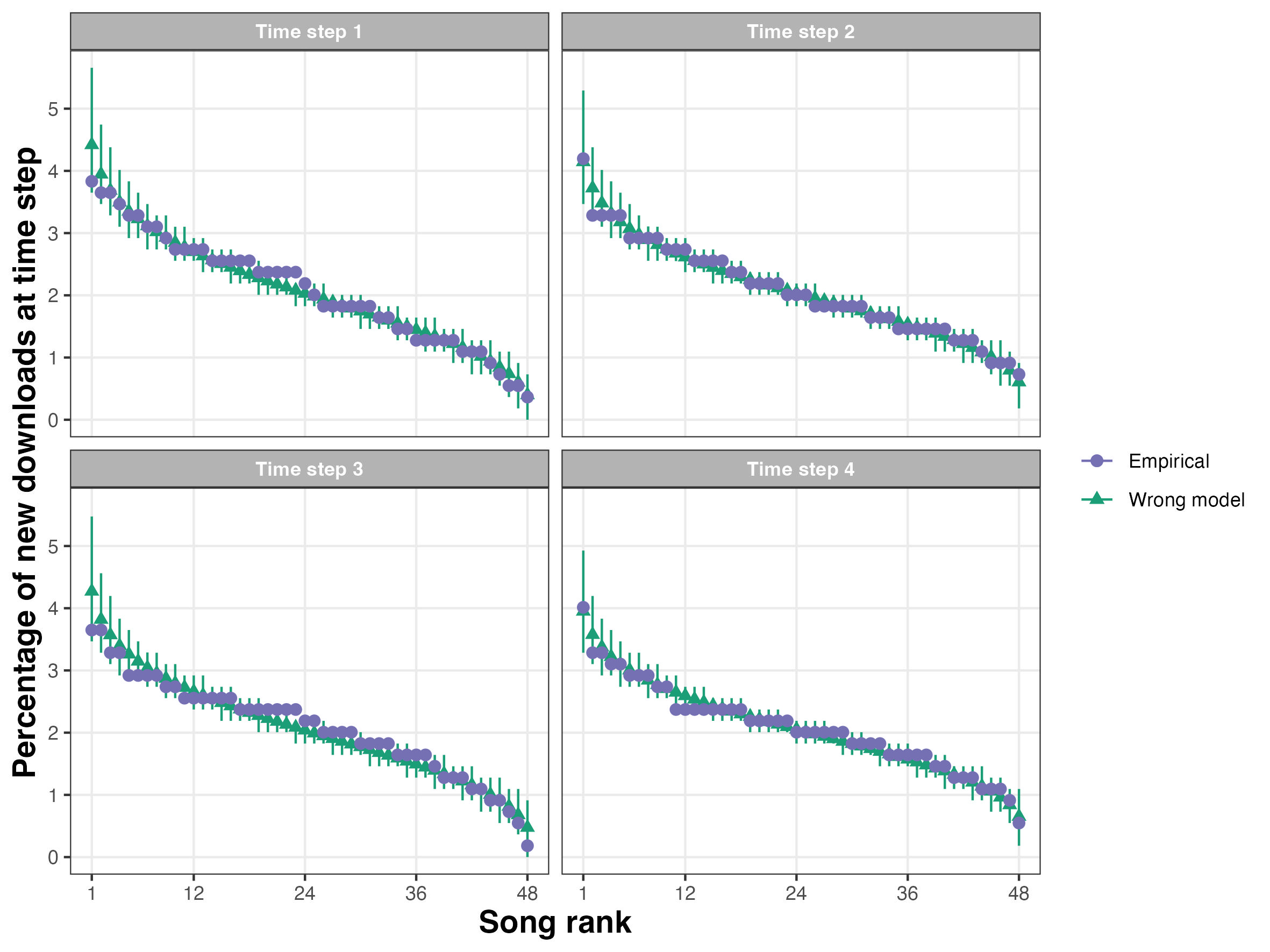}
    \caption{Same as \cref{fig:musicLab1}, for version 2 of the MusicLab experiment. Here \textit{f}=0.310.}
    \label{fig:musicLab2}
\end{figure*}

\section*{Discussion}
The results are cause for a reassessment of the large body of empirical evidence on the origins of human success accumulated across literatures and contexts. They call into question previous conclusions drawn about the existence of cumulative advantage processes based on longitudinal analysis of individual achievement records \cite{diprete2006cumulative}. The present paper implies that these records equally fit a model that lacks any notion of path dependence. Vice versa, the data used in studies that fit models that base success on pre-existing individual qualities equally fit predictions from a model of path dependence in which all individuals are ex ante equally talented \cite{huber1998cumulative, sinatra2016quantifying}.

How forward? The principal way out of the fundamental confoundedness of path dependence and talent differences in longitudinal data is (quasi-)experimental methods. Multiple world experiments can reveal path dependence by demonstrating varying outcomes under identical talent conditions \cite{salganik2006experimental, frey2021social, centola2015spontaneous}. Alternatively, regression discontinuity designs, natural experiments, and matched case control studies can be used to test the success-breeds-success hypothesis by examining effects of random success on subsequent success \cite{bol2018matthew, azoulay2014matthew, sorensen2007bestseller, wang2019early, chan2014academic, borjas2015prizes, farys2021matthew, farys2017matched}.

In the absence of such (quasi-)experimental control, studies working with temporal success data should be cognizant of the existence of twin models with and without talent or path dependence. Past methods that have been explicitly developed for disentangling path dependence from talent differences as distinct sources of success necessarily make auxiliary assumptions in order to attribute variance to either source \cite{allison1980estimation, holden1986contagiousness, huber2002new, box2006repeated}. When these models are used, these auxiliary assumptions should be made explicit and grounded in theory. In some cases one model may be preferred over the other in terms of parsimony. For example, the $Q$-model appears to be more parsimonious over its path dependent twin that we identified here, because the latter lacks a mechanism to justify the decreasing conditional variance in \cref{eqVariance1}, while in the $Q$-model this is a natural consequence of the fact that uncertainty about one's talent decreases as they produce more publications. Vice versa, if this conditional variance were empirically found to be constant instead, then the corresponding path dependent model would be more parsimonious than its talent equivalent.

\bibliography{nonidentifiability}

\appendix
\section*{Appendix A - Examples of purely path dependent and talent models}\label{a:appendix1}

Here we show how some standard models of cumulative advantage can be described within our framework, i.e. written as purely path dependent models according to our definition. We also give an example of a model from the literature (the random-typing monkey) that can be expressed as a talent model according to our definition. In what follows, by ${\mathcal U}(0,1)$ we denote the uniform probability distribution on the interval $[0,1]$.

\subsection*{P\'olya urn}

In the standard P\'olya urn, the probability that a red ball is added at step $n$ is equal to $\frac{n_R}{n-1+K}$, where $n_R$ is the number of red balls currently in the urn (at time $n-1$) and $K$ is the initial number of balls in the urn. Letting $Y_n=1$ denote the fact that a red ball is drawn at step $n$ (and $Y_n=0$ otherwise), then $n_R$ is equal to $\sumdu{k=1}{n-1}Y_k+K_R$, where $K_R$ is the initial number of red balls. Therefore, the probability of the event $Y_n=1$ is equal to $\frac{\sumdu{k=1}{n-1}Y_k+K_R}{n-1+K}$.

We now define a path dependent process, according to the definition given in the current manuscript, which is equivalent to the P\'olya urn just described. Let $\{X_n\}$ be an i.i.d. sequence of $\mathcal{U}(0,1)$ random variables and define $Y_n=g_n(Y_1,\ldots ,Y_{n-1},X_n)$, where
\begin{equation}
    g_n(Y_1,\ldots ,Y_{n-1},X_n)=\left\{\begin{matrix}
        1, & \text{if } X_n<\frac{\sumdu{k=1}{n-1}Y_k+K_R}{n-1+K}\\
        0, & \text{otherwise.}
    \end{matrix}\right.
\end{equation}

Clearly, conditioned on $Y_1,\ldots ,Y_{n-1}$, the probability of the event $Y_n=1$ is equal to $\frac{\sumdu{k=1}{n-1}Y_k+K_R}{n-1+K}$, exactly as required. We thus see that the P\'olya urn model fits our definition of a path dependent process.

\subsection*{Preferential attachment}

A preferential attachment model is similar to a P\'olya urn, except that balls of new colors may enter the urn. For example, Simon \cite{simon1955class} describes the generation of a text by adding words one by one, with the probability of adding a given word being proportional to the occurrences of that word in the text, except that there is a probability $\alpha\in (0,1)$ that an entirely new word appears. If we focus on a particular word that is already present in the text,\footnote{In this model, the way that newly appearing words are chosen is not specified. One can say that words do not exist before their first appearance. Thus, specifying the probabilities of the first appearance of a specific word in this model is impossible without further assumptions.} the probability of a new occurrence of this word at step $n$ is then $(1-\alpha)\cdot \frac{n_R}{n-1+K}$, where $n_R$ is the number of occurrences of that word so far, and $K$ is the number of words in the text at time $0$. This model can be expressed in terms of our definition of a path dependent process in a way similar to a P\'olya urn.

Barab\'asi and Albert \cite{barabasi1999emergence} describe the growth of a network in which at each step a new node is added together with a link, which connects it to an existing node. The connecting node is chosen randomly, with probability proportional to its current degree, i.e. the number of existing links to that node. Assuming that the network at time $t=0$ contains $m_0$ links, then at time $t$ it has $t+m_0$ links, and because each link contributes to the degree of two nodes, the sum of the degrees is $2(t+m_0)$. A node that enters at time $t'$ will therefore have zero probability of obtaining links for $t<t'$, probability equal to $1$ at time $t=t'$, and probability $\frac{X_t}{2(t+m_0)}$ of obtaining a link at any $t>t'$, where $X_t$ is its degree. Again, this model can be cast in terms of our definition of a path dependent model in a way similar to a P\'olya urn.

We note however that the resulting probability distribution of $X_t$ is not the distribution of the degree of a randomly chosen node at time $t$, because nodes enter at different times. What the distribution of $X_t$ tells us is the probability that the degree of the node entering at the pre-specified time $t'$ will take any given value at time $t$.

\subsection*{Price's urn model}

Price \cite{price1976general} defines the following urn model: An urn with 1 black and 1 red ball is sampled repeatedly, with red balls interpreted as successes. Whenever a red ball is drawn, it is returned together with another red ball. The first time a black ball is drawn, we stop drawing balls altogether. Interpreting $Y_n=1$ as a success (red ball drawn) and $Y_n=0$ as a failure (black ball drawn or no draw), we may define an equivalent generative model as follows: $X_n$ are i.i.d. $\mathcal{U}(0,1)$ random variables and $Y_n=g_n(Y_1,\ldots ,Y_{n-1},X_n)$, where
\begin{equation}
    g_n(Y_1,\ldots ,Y_{n-1},X_n)=\left\{\begin{matrix}
        1, & \text{if } Y_1=\ldots =Y_{n-1}=1\\
        & \text{ and }X_n<\frac{n}{n+1},\\
        0, & \text{otherwise.}
    \end{matrix}\right.
\end{equation}

This is again a special case of a path dependent model according to our definition.

\subsection*{Gibrat's Law of proportional effect}

Gibrat \cite{gibrat1931inegalits} suggested that the size of firms in terms of number of employees, as well as the size of cities in terms of number of inhabitants, grows on average proportionally to their current size. Specifically, if $Y_{n-1}$ is the size of a firm in year $n-1$, then the size in year $n$ is given by $Y_{n}=Y_{n-1}\cdot (1+X_n)$ where $\{X_n\}$ is an i.i.d. sequence of random variables. Thus, this is a special case of a path dependent process according to our definition, where $g_n(Y_1,\ldots ,Y_{n-1},X_n)=Y_{n-1}\cdot (1+X_n)$.

The above can be generalized to the case that the $X_n$'s are not i.i.d., even though our definition of a path dependent process requires an i.i.d. sequence $X_n'$. To do so, we may choose $X_n'$ to be i.i.d. uniformly distributed in $[0,1]$ and define $g_n(Y_1,\ldots ,Y_n,X_n')=Y_n\cdot (1+F_{n}^{-1}(X_n'))$, where $F_n^{-1}$ is the generalized inverse of the cumulative distribution function of $X_n$. By a standard result in statistics (see inverse transform sampling), the random variable $F_{n}^{-1}(X_n')$ has the same distribution as $X_n$.

\subsection*{Random-typing monkey}

Miller \cite{miller1957some} showed that a power law distribution for the number of occurrences of different words in a text can be obtained if we simply assume that a monkey presses keys at random. Specifically, we assume that each key is pressed with the same probability, except perhaps from the space key, which is pressed with some fixed probability $\alpha \in (0,1)$. If there are $N$ different keys other than the space key, then the probability of forming a specific word of length $k$ is $\left(\frac{1-\alpha}{N}\right)^{k}\cdot \alpha$. Because this probability does not depend on previous occurrences of the word, this is a talent model, specifically a stationary one.

To cast this in the form of a talent model according to the definition we have given, we let $\{X_n\}$ be i.i.d. random variables distributed as $\mathcal{U}(0,1)$ and define $T$ to be the length of the word.\footnote{Here larger $T$ makes success less likely, so $T$ is negative talent.} The function $f$ is given by
\begin{equation}
    f(T,X_n)=\left\{\begin{matrix}
        1, & \text{if } X_n<\left(\frac{1-\alpha }{N}\right)^T\cdot \alpha \\
        0, & \text{otherwise.}
    \end{matrix}\right.
\end{equation}
We get that conditioned on word length $T=k$, the probability of appearance of the word ($f(T,X_n)=1$) at step $n$ is $\left(\frac{1-\alpha }{N}\right)^k\cdot \alpha $, as required.

The above assumes that $T$ is given, which means that we follow a predefined word. To be consistent with our definition, $T$ must be drawn from a given probability distribution, which here it would be equivalent to choosing stochastically a word to track (observe its occurrences). For example, if the word we tracked was chosen uniformly randomly from among all possible words of length at most $K$, then $T$ would equal $k\in \{1,\ldots ,K\}$ with probability $\frac{N^k}{\sumdu{l=1}{K}N^l}$, because there are exactly $N^k$ words of length $k$.

\section*{Appendix B - Proofs of theorems}

\begin{proof}[Proof of \cref{theoremTalentToPathDependent}]
Let $\{X_n\}$ be i.i.d. random variables distributed as $\mathcal{U}(0,1)$. Suppose that $g_1,\ldots ,g_{n-1}$ have already been defined such that $(Y'_1,\ldots ,Y'_{n-1})$ has the same distribution as $(Y_1,\ldots ,Y_{n-1})$, where $Y'_k=g_k(Y'_1,\ldots ,Y'_{k-1},X_k)$ for $k\in \{1,\ldots ,n-1\}$. By the transfer theorem \cite[Theorem 6.10]{kallenberg2002foundations}, there exists some measurable function $g_n$, such that even $(Y'_1,\ldots ,Y'_n)$ has the same distribution as $(Y_1,\ldots ,Y_n)$, where $Y'_n=g_n(Y'_1,\ldots ,Y'_{n-1},X_n)$.
\end{proof}

\begin{proof}[Proof of \cref{propositionQmodel}]
    The variables $Y_1,\ldots ,Y_n$ can be thought of as noisy observations of the variable $T$, which has prior distribution ${\mathcal{N}}(\mu_T,\sigma^2_T)$, and the noise is distributed according to ${\mathcal{N}}(0,\sigma^2_X)$. By a standard result in Bayesian estimation (see for example \cite[Ch. 4, Example 2.2]{lehmann2006theory}), conditioned on $Y_1,\ldots ,Y_n$, the posterior distribution of $T$ is normal with mean
    $$
    \frac{\frac{\sigma^2_X}{\sigma_T^2}\cdot \mu_T+\sumdu{k=1}{n}Y_k}{n+\frac{\sigma^2_X}{\sigma_T^2}}
    $$
    and variance $\frac{\sigma^2_X}{n+\frac{\sigma^2_X}{\sigma_T^2}}$. Therefore, conditioned on $Y_1,\ldots, Y_n$, the variable $Y_{n+1}=T+X_{n+1}$ is distributed normally with the same mean and with variance
    $$
    \frac{\sigma^2_X}{n+\frac{\sigma^2_X}{\sigma_T^2}}+\sigma ^2_X=\sigma ^2_X\cdot \left(1+\frac{1}{n+\frac{\sigma^2_X}{\sigma_T^2}}\right).
    $$
    This is the same as the conditional distribution of $Y'_{n+1}$, conditioned on $Y'_1,\ldots, Y'_n$, if we set  $a=\mu_T$, $b=\sigma ^2_X$, and $c=\frac{\sigma^2_X}{\sigma_T^2}$ (by \cref{eqMean1,eqVariance1}).

    Since the conditional distribution of $Y_{n+1}$ given $Y_1,\ldots ,Y_n$ is the same as the conditional distribution of $Y'_{n+1}$ given $Y'_1,\ldots, Y'_n$ for all $n\in\N$, the joint distributions of $Y_1,Y_2,\ldots $ and $Y'_1,Y'_2,\ldots $ are identical.
\end{proof}

\begin{proof}[Proof of \cref{theoremPathDependentToTalentExchangeable}]
    For any talent model $Y_n'=f(T,X_n)$, conditioned on $T$, the sequence $\{Y_n'\}$ is i.i.d., because $\{X_n\}$ is. Therefore, $\{Y_n'\}$ is a mixture of i.i.d. sequences, with mixing distribution that of $T$. By the converse of De Finetti's theorem, $\{Y_n'\}$ is exchangeable.
    
    For the converse, suppose that $\{Y_n\}$ is exchangeable. By De Finetti's theorem, it is a mixture of i.i.d. sequences, i.e. there is some random variable $T$, independent of $\{X_n\}$, such that conditioned on $T$, the $Y_n$'s are i.i.d. Let $\{X_n\}$ be $\mathcal{U}(0,1)$ i.i.d. By the randomization lemma \cite[Lemma 3.22]{kallenberg2002foundations}, there exists some measurable function $f$ such that $f(T,X_1)$ has the same distribution as $Y_1$. Define $Y_n'=f(T,X_n)$. Since $\{X_n\}$ is i.i.d. and independent of $T$, we get that $Y_n'=f(T,X_n)\sim f(T,X_1)\sim Y_1\sim Y_n$. Therefore, conditioned on $T$, both sequences $\{Y_n\}$ and $\{Y_n'\}$ are i.i.d. and with the same distribution. We conclude that $\{Y_n\}\sim \{Y_n'\}$ even unconditionally.
\end{proof}

\begin{proof}[Proof of \cref{theoremPathDependentToTalent}]
    Let $T$ be $\mathcal{U}(0,1)$. By the randomization lemma \cite[Lemma 3.22]{kallenberg2002foundations}, there exists some measurable function $f:T\mapsto (Y_1',Y_2',\ldots )$ such that $\{Y_n'\}$ has the same distribution as $\{Y_n\}$. Denoting the $n$-th component of $f$ by $f_n$, we have that $Y'_n=f_n(T)$.
\end{proof}

\begin{proof}[Proof of \cref{propAllisonTalentEquivalent}]
    Define a time-rescaled version of the contagious Poisson process, using the transformation $s(t)={\frac{1}{\beta }\cdot \ln (\beta t+1)}$, i.e. $Y_t=X_{s(t)}$. The inverse of this transformation is $g(t)=\frac{1}{\beta }\cdot (e^{\beta t}-1)$. The process $Y_t$ is again a Markov counting process with rate
    \begin{equation}
        s'(t)\cdot \left(\alpha +\beta X_{s(t)}\right)=\frac{1}{\beta t+1}\cdot (\alpha +\beta Y_t)=\frac{\frac{\alpha}{\beta}+Y_t}{\beta t+1}\cdot \beta.
    \end{equation}
    By \cite{cane1974concept}, this is equivalent to the counting process associated with a mixed Poisson process with (homogeneous) rate $\lambda\beta $, where $\lambda \sim \text{Gamma}\left(\frac{\alpha}{\beta },1\right)$, or equivalently, with rate distributed according to $\text{Gamma}\left(\frac{\alpha}{\beta },\beta\right)$, by the scaling property of the Gamma distribution. Therefore, the process $X_t=Y_{g(t)}$ is (the counting process of) a non-homogeneous mixed Poisson process with rate $\lambda \cdot g'(t)=\lambda \cdot e^{\beta t}$, where $\lambda $ is distributed according to $\text{Gamma}\left(\frac{\alpha}{\beta },\beta \right)$.
\end{proof}

\section*{Appendix C - Simulation methods}

\subsubsection*{Notes on the P\'olya urn simulations}

The value of the parameter $f$ was chosen so that to minimize the sum of absolute errors between the empirical data and the mean download proportions for each rank in 300 simulations. Specifically, we varied $f$ in the interval $0.1-0.6$ with steps of $0.05$.

The simulations in the two versions of the experiment differ in the number of participants and their respective number of downloads. This is why the mean and the error bars displayed in \cref{fig:musicLab1,fig:musicLab2} are not exactly the same. The code is available online at \url{https://github.com/Lsage/Cumulative-advantage/}.

\end{document}